\documentclass[12pt]{amsart}
\usepackage[utf8]{inputenc}
\usepackage[T1]{fontenc}
\usepackage{amsfonts}
\usepackage[leqno]{amsmath}
\usepackage{accents}
\usepackage{amssymb, comment, float}
\usepackage[dvipsnames]{xcolor}
\usepackage[foot]{amsaddr}
\usepackage[shortlabels]{enumitem}
\usepackage{mathdots}
\usepackage[a4paper, footskip=0.5in,  headheight = 0.5in, top=1.25in, bottom=1.25in,  right=1in,  left=1in]{geometry}
\usepackage{hyperref}
\hypersetup{
colorlinks,
linkcolor=black,
urlcolor=Blue,
citecolor=black,}
\usepackage{cleveref}
\usepackage[center]{caption}
\usepackage{eufrak}
\usepackage{marvosym}     
\usepackage{latexsym}
\usepackage[nodayofweek]{datetime}
\usepackage{tikz}
\usepackage{subfig}
\usepackage{thm-restate}
\usepackage{verbatim}

\DeclareFontFamily{U}{mathx}{\hyphenchar\font45}
\DeclareFontShape{U}{mathx}{m}{n}{
      <5> <6> <7> <8> <9> <10>
      <10.95> <12> <14.4> <17.28> <20.74> <24.88>
      mathx10
      }{}
\DeclareSymbolFont{mathx}{U}{mathx}{m}{n}
\DeclareFontSubstitution{U}{mathx}{m}{n}
\DeclareMathAccent{\widecheck}{0}{mathx}{"71}
\DeclareMathAccent{\wideparen}{0}{mathx}{"75}

\newtheorem{statement}{statement}[section]
\newtheorem{theorem}[statement]{Theorem}
\newtheorem{lemma}[statement]{Lemma}

\newtheorem{corollary}[statement]{Corollary}

\DeclareMathOperator{\tw}{tw}
\DeclareMathOperator{\ta}{{\text{tree}-\alpha}}

\newcommand{\mca}{\mathcal}

\newcommand{\poi}{\mathbb{N}}

\newcounter{tbox}
\newcommand{\sta}[1]{\vspace*{0.5cm}\refstepcounter{tbox}\noindent{\parbox{\textwidth}{(\thetbox) \emph{#1}}}\vspace*{0.3cm}}
\newcommand{\mylongtitle}[1]{%
  \ifodd\value{page}%
    \protect\parbox{0.97\linewidth}{#1}\hfill%
  \else%
    \hfill\protect\parbox{0.97\linewidth}{#1}%
  \fi%
}

\title[Induced subgraphs and tree decompositions XIX.]{Induced subgraphs and tree decompositions\\
XIX. thetas and forests}

\author{Maria Chudnovsky$^{\ast \dagger}$}
\author{Julien Codsi$^{\ast \ddagger}$}
\author{Sepehr Hajebi$^{\mathsection}$}
\author{Sophie Spirkl$^{\mathsection \parallel}$}

\thanks{$^{\ast}$ Princeton University, Princeton, NJ, USA}
\thanks{$^{\mathsection}$ Department of Combinatorics and Optimization, University of Waterloo, Waterloo, Ontario, Canada}
\thanks{$^{\dagger}$ Supported by NSF Grant DMS-2348219, NSF Grant CCF-2505100, AFOSR grant 
FA9550-22-1-0083 and a Guggenheim Fellowship.}
\thanks{$^{\ddagger}$  Supported by 
NSF Grant DMS-2348219, 
AFOSR grant  FA9550-22-1-0083  and  the Fonds de recherche du Qu\'{e}bec Grant 321124.} 
\thanks{$^{\parallel}$ We acknowledge the support of the Natural Sciences and Engineering Research Council of Canada (NSERC), [funding reference number RGPIN-2020-03912].
Cette recherche a \'et\'e financ\'ee par le Conseil de recherches en sciences naturelles et en g\'enie du Canada (CRSNG), [num\'ero de r\'ef\'erence RGPIN-2020-03912]. This project was funded in part by the Government of Ontario. This research was conducted while Spirkl was an Alfred P. Sloan Fellow.}
\date {\today}
\begin{document}
\maketitle

\begin{abstract}
Let $H$ be a graph and let $\mca{C}$ be a hereditary class of theta-free graphs such that $H\notin \mca{C}$. We prove that if
\begin{itemize}
 \item[(a)] $H$ is a forest; and
\item[(b)] $\mca{C}$ excludes the line graphs of all subdivisions of some wall;
\end{itemize}
then the treewidth of every graph in $\mca{C}$ is at most a polynomial function of its clique number. This is best possible in that both (a) and (b) are necessary for the existence of {\sl any} function with the above property.
 \end{abstract}
 
\section{Introduction}
The set of positive integers is denoted by $\poi$. For integers $k,k'$, we write $\{k,\ldots, k'\}$ for the set of integers no smaller than $k$ and no larger than $k'$. Graphs in this paper have finite vertex sets, no loops, and no parallel edges. Given a graph $G=(V(G),E(G))$ and $X\subseteq V(G)$, we use $X$ to denote both the set $X$ of vertices and the induced subgraph $G[X]$ of $G$ with vertex set $X$. A \textit{clique} in $G$ is a set of pairwise adjacent vertices, and a \textit{stable set} is a set of pairwise nonadjacent vertices. The \textit{clique number} of $G$, denoted $\omega(G)$, is the maximum cardinality of a clique in $G$. For a graph $H$, we say that $G$ is \textit{$H$-free} if $G$ has no induced subgraph isomorphic to $H$. For two graphs $H_1, H_2$, we say that $G$ is \textit{$(H_1, H_2)$-free} if $G$ is both $H_1$-free and $H_2$-free. The \textit{treewidth} of a graph $G$ is denoted by $\tw(G)$ (see \cite{diestel} for the definition of treewidth and other standard terminology).

\begin{figure}
    \centering
    \includegraphics[scale=0.8]{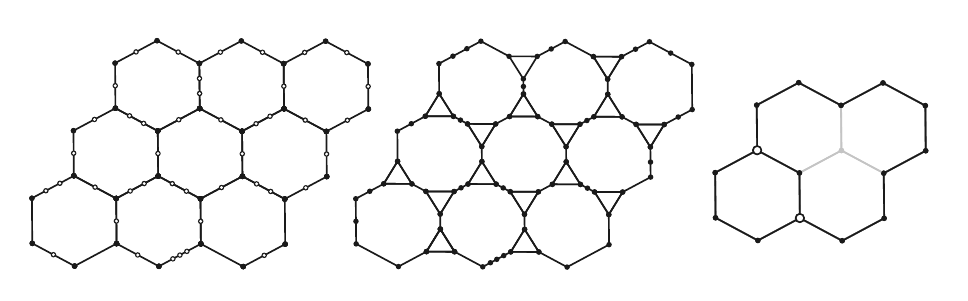}
    \caption{A subdivision of $W_{4\times 4}$ (left), its line graph (middle), and a theta in $W_{3\times 3}$ (right).}
    \label{fig:basic}
\end{figure}

This series of papers studies the following question: What induced subgraphs must be excluded to guarantee bounded treewidth? In principle, this is motivated by the grid theorem of Robertson and Seymour, \Cref{thm:RSPW} below, which gives a complete answer to the analogous question for minors (and subgraphs). For every $t\in \poi$, we denote by $W_{t\times t}$ the $t$-by-$t$ hexagonal grid, also known as the \textit{$t$-by-$t$ wall}; see \Cref{fig:basic}. It is well known \cite{diestel} that every subdivision of $W_{t\times t}$ has treewidth $t$.

\begin{theorem}[Robertson and Seymour \cite{GMV}]\label{thm:RSPW}
For every $t\in \poi$, graphs with no minor (or equivalently, no subgraph) isomorphic to any subdivision of $W_{t\times t}$ have bounded treewidth.
\end{theorem}

In contrast, excluding subdivided walls as induced subgraphs does not necessarily bound the treewidth, nor does excluding a much simpler family of graphs: thetas. A \textit{theta} is a graph that is isomorphic to a subdivision of the complete bipartite graph $K_{2,3}$. A \textit{theta in} a graph $G$ is an induced subgraph of $G$ that is a theta, and we say that $G$ is \textit{theta-free} if there is no theta in $G$. There is a theta in every subdivision of $W_{r\times r}$ for all $r\geq 3$ (see \Cref{fig:basic}), and thus theta-freeness is a much stronger restriction than excluding all subdivisions of some fixed wall as induced subgraphs. That said, there are theta-free graphs of arbitrarily large treewidth: Complete graphs and the line graphs of subdivided walls (see \Cref{fig:basic}) are theta-free, and it is well known \cite{diestel} that for each $r\in \poi$, both $K_{r+1}$ and the line graphs of all subdivisions of $W_{r\times r}$ have treewidth $r$.

What about theta-free graphs that also exclude both large complete graphs and the line graphs of subdivisions of large walls? The easiest nontrivial case is when the graphs in question are triangle-free (and thus exclude the line graphs of subdivided walls, since those contain triangles). Sintiari and Trotignon showed that even under this assumption, and much stronger ones, the treewidth remains unbounded (recall that the \textit{girth} of a graph $G$ is the length of the shortest cycle in $G$):

\begin{theorem}[Sintiari and Trotignon \cite{layeredwheels}]\label{thm:layeredwheel}
    For all $g,w\in \poi$, there is a theta-free graph of girth larger than $g$ and treewidth larger than $w$.
\end{theorem}

The graphs used to prove \Cref{thm:layeredwheel} are an example of a rich class of constructions known as ``layered wheels'' \cite{newlayeredwheels, CT, layeredwheels}. Layered wheels are often rather involved, and in a sense, layered-wheel-like obstructions to bounded treewidth are the last remaining hurdle to proving a full analog of \Cref{thm:RSPW} for induced subgraphs (where by ``layered-wheel-like obstructions'' we mean, roughly, those which do contain large complete minors but do not contain large complete bipartite \textit{induced} minors; see \cite{tw16} for more details).

Therefore, we are very much interested in understanding the structural properties of layered wheels. For example, what induced subgraphs are unavoidable in layered wheels of large enough treewidth? One may observe \cite{layeredwheels} that every forest $H$ is contained as an induced subgraph in all theta-free layered wheels (meaning those from \Cref{thm:layeredwheel}) of sufficiently large treewidth. On the other hand, every graph $H$ with the latter property must be a forest because theta-free layered wheels can have arbitrarily large girth. More generally, the following is an immediate corollary of \Cref{thm:layeredwheel}. 

\begin{corollary}\label{cor:onlyif}
    Let $H$ be a graph such that for all $r\in \poi$, every theta-free $H$-free graph with bounded clique number and no induced subgraph isomorphic to the line graph of any
 subdivision of $W_{r\times r}$ has bounded treewidth. Then $H$ is a forest.
\end{corollary}

Our main result is the converse to \Cref{cor:onlyif}, with the treewidth bound polynomial in the clique number.

\begin{restatable}{theorem}{maintheorem}\label{thm:main}
For all $r\in \poi$ and every forest $H$, there is a constant $d_{\ref{thm:main}}=d_{\ref{thm:main}}(H,r)\in \poi$ such that for every $t\in \poi$, every theta-free $(H, K_t)$-free graph $G$ with no induced subgraph isomorphic to the line graph of any subdivision of $W_{r\times r}$ satisfies $\tw(G)\leq t^{d_{\ref{thm:main}}}$.
\end{restatable}

In other words, \Cref{thm:main} suggests that theta-free layered wheels, although globally intricate, are locally canonical, in that they represent the local structure of \textit{all} theta-free graphs of large treewidth (that avoid complete graphs and the line graphs of subdivided walls). 

From \Cref{thm:main}, we deduce the following (note that ``\ref{cor:main_c} implies \ref{cor:main_b}'' and ``\ref{cor:main_b} implies \ref{cor:main_a}'' are trivial, and ``\ref{cor:main_a} implies \ref{cor:main_c}'' follows from \Cref{thm:main}).

\begin{corollary}\label{cor:main}
    Let $H$ be a forest and let $\mca{C}$ be a hereditary class of theta-free graphs with $H\notin \mca{C}$. Then the following are equivalent.
  \begin{enumerate}[{\rm(a)}]
  \item\label{cor:main_a} There exists $r\in \poi$ such that for every subdivision $W$ of $W_{r\times r}$, the line graph of $W$ does not belong to $\mca{C}$. 
  \item\label{cor:main_b} There is a function $f:\poi\to\poi$ such that $\tw(G)\leq f(\omega(G))$ for every $G\in \mca{C}$. 
   \item\label{cor:main_c} There is a polynomial function $f:\poi\to\poi$ such that $\tw(G)\leq f(\omega(G))$ for every $G\in \mca{C}$. 
  \end{enumerate}
\end{corollary}

\Cref{thm:main} is also a substantial generalization of the main result of an earlier paper in this series \cite{tw8}. A \textit{prism} is the line graph of a theta, and a graph $G$ is \textit{prism-free} if no induced subgraph of $G$ is a prism. For instance, if $r\geq 3$, then the line graphs of subdivisions of $W_{r\times r}$ are not prism-free.

\begin{theorem}[Abrishami, Alecu, Chudnovsky, Hajebi, Spirkl \cite{tw8}]\label{thm:tw8}
For all $t\in \poi$ and every forest $H$, there is a constant $c_{\ref{thm:tw8}}=c_{\ref{thm:tw8}}(H,t)$ such that every graph $G$ that is theta-free, prism-free, and $(H,K_t)$-free satisfies $\tw(G)\leq c_{\ref{thm:tw8}}$.
\end{theorem}

As another corollary, from \Cref{thm:main} combined with the main result of \cite{logDMS} (see Theorem~1.1 in \cite{logDMS}), it follows that for all $r\in \poi$ and every forest $H$,  theta-free $H$-free graphs with no induced subgraph isomorphic to the line graph of a subdivision of $W_{r\times r}$ have poly-logarithmic ``tree independence number.'' The \textit{tree independence number} of a graph $G$, denoted $\ta(G)$, is the minimum $s\in \poi$ for which $G$ admits a tree decomposition where each bag intersects every stable set of $G$ in at most $s$ vertices.

\begin{corollary}\label{cor:treealpha}
    For all $r\in \poi$ and every forest $H$, there are constants $c_{\ref{thm:main}}=c_{\ref{thm:main}}(H,r)\in \poi$ and $d_{\ref{thm:main}}=d_{\ref{thm:main}}(H,r)\in \poi$ such that every theta-free $H$-free graph $G$ with $|V(G)|>1$ and no induced subgraph isomorphic to the line graph of any subdivision of $W_{r\times r}$ satisfies $\ta(G)\leq c_{\ref{thm:main}}\log^{d_{\ref{thm:main}}} (|V(G)|)$.
\end{corollary}

\Cref{cor:treealpha}, in turn, implies that the \textsc{Maximum Weight Independent Set} problem (and many other problems that are {\sf NP}-hard in general) can be solved in quasi-polynomial time in every class $\mca{C}$ as above (see \cite{ti4} and \cite{logDMS} for more details).

\section{Reducing to low separability}

The key tool in the proof of our main result is \Cref{thm:mainseparation} below. In this section, we show that \Cref{thm:main} follows from \Cref{thm:mainseparation}. The remainder of the paper will then be devoted to the proof of \Cref{thm:mainseparation}, which we will complete in \Cref{sec:final}.

We begin with some definitions. Let $G$ be a graph. A \textit{path in $G$} is an induced subgraph $P$ of $G$ that is a path. The vertices of $P$ with degree at most one in $P$ are called the \textit{ends of $P$}, and the set of all other vertices of $P$ is called the \textit{interior of $P$}, denoted $P^*$. The \textit{length} of a path in $G$ is its number of edges. If $P$ is a path in $G$ with ends $x,y\in V(G)$, then we say that $P$ is a path in $G$ \textit{from $x$ to $y$}. For $\lambda\in \poi$, we say that $G$ is \textit{$\lambda$-separable} if for every two distinct and nonadjacent vertices $x,y\in V(G)$, there is no set of $\lambda$ pairwise internally disjoint paths in $G$ from $x$ to $y$. 

\begin{restatable}{theorem}{mainseparation}\label{thm:mainseparation}
    For all $h\in \poi$, there is a constant $d_{\ref{thm:mainseparation}}=d_{\ref{thm:mainseparation}}(h)\in \poi$ such that for every forest $H$ on at most $h$ vertices and every $t\in \poi$, every theta-free $(H, K_{t})$-free graph is $t^{d_{\ref{thm:mainseparation}}}$-separable. 
\end{restatable}

In order to derive \Cref{thm:main} from \Cref{thm:mainseparation}, we use several results from the literature. First, we need a result from \cite{polypw}. An \textit{induced minor} of a graph $G$ is a graph obtained from $G$ by deleting vertices and contracting edges, and removing the loops and the parallel edges arising in this process.

\begin{theorem}[Hajebi; Theorem~3.2 in \cite{polypw} for $\kappa=2$]\label{thm:mainpolyblock}
    For every planar graph $H$ and every $\sigma\in \poi$, there is a constant $d_{\ref{thm:mainpolyblock}}=d_{\ref{thm:mainpolyblock}}(H,\sigma)\in \poi$ such that for all $\lambda,t\in \poi$, if $G$ is a $\lambda$-separable $K_t$-free graph with no induced minor isomorphic to $H$ or $K_{\sigma,\sigma}$, then $\tw(G)\leq (2\lambda t)^{d_{\ref{thm:mainpolyblock}}}$.
\end{theorem}

Next, we need the following from \cite{tw16}. For $s,l\in \poi$,  an \textit{$(s,l)$-constellation} is a graph $C$ in which there is a stable set $S$ with $|S|=s$ such that $C\setminus S$ has exactly $l$ components, every component of $C\setminus S$ is a path, and each vertex $x\in S$ has at least one neighbor in $C$ in each component of $C\setminus S$.  Given a graph $G$, by an \textit{$(s,l)$-constellation in $G$} we mean an induced subgraph of $G$ which is an $(s,l)$-constellation.

\begin{theorem}[Chudnovsky, Hajebi, Spirkl; Theorem~1.2 in \cite{tw16}]\label{thm:compbip}
   For all $l,r,s\in \poi$, there is a constant $c_{\ref{thm:compbip}}=c_{\ref{thm:compbip}}(l,r,s)\in \poi$ such that for every graph $G$ with an induced minor isomorphic to $K_{c_{\ref{thm:compbip}},c_{\ref{thm:compbip}}}$, one of the following holds: 
\begin{enumerate}[{\rm (a)}]
  \item\label{thm:compbip_a} $G$ has an induced minor isomorphic to $W_{r\times r}$.
  \item\label{thm:compbip_b} There is a $(s,l)$-constellation in $G$.
\end{enumerate} 
\end{theorem}

We also need a well-known result from \cite{aboulker}:

\begin{lemma}[Aboulker, Adler, Kim, Sintiari, Trotignon; Corollary of Lemma 3.6 in \cite{aboulker}] \label{lem:minor-to-ind}
    For every $r \in \poi$, there is a constant $c_{\ref{lem:minor-to-ind}} = c_{\ref{lem:minor-to-ind}}(r)$ such that if $G$ is a graph with an induced minor isomorphic to $W_{c_{\ref{lem:minor-to-ind}} \times c_{\ref{lem:minor-to-ind}}}$, then $G$ has an induced subgraph isomorphic to either a subdivision of $W_{r\times r}$ or the line graph of a subdivision of $W_{r\times r}$.
\end{lemma}

We are now ready to derive \Cref{thm:main} from \Cref{thm:mainseparation}:

\begin{proof}[Proof of \Cref{thm:main} {\sl(Assuming \Cref{thm:mainseparation})}]

Let 
$$\rho=c_{\ref{lem:minor-to-ind}}(\max\{r,3\});\quad 
\sigma=c_{\ref{thm:compbip}}(3,\rho,2);\quad 
d=d_{\ref{thm:mainseparation}}(|V(H)|);\quad d'=d_{\ref{thm:mainpolyblock}}(W_{\rho\times \rho},\sigma).$$

We will show that 
$$d_{\ref{thm:main}}=d_{\ref{thm:main}}(H,r)=(d+2)d'$$
works. Let $t\in \poi$ and let $G$ be a theta-free $(H, K_t)$-free graph $G$ with no induced subgraph isomorphic to the line graph of any subdivision of $W_{r\times r}$. By \Cref{thm:mainseparation} $G$ is $t^{d}$-separable. We further claim that:

\sta{\label{st:nowall} $G$ has no induced minor isomorphic to $W_{\rho\times \rho}$.}

Note that $G$ has no induced subgraph isomorphic to any subdivision of $W_{3\times 3}$ since $G$ is theta-free, and recall that $G$ has no induced subgraph isomorphic to the line graph of a subdivision of $W_{r\times r}$. Thus, $G$ has no induced subgraph isomorphic to a subdivision of $W_{\max\{r,3\}\times \max\{r,3\}}$ or the line graph of such a subdivision of $W_{\max\{r,3\}\times \max\{r,3\}}$. It follows from  \Cref{lem:minor-to-ind} and the choice of $\rho$ that $G$ has no induced minor isomorphic to $W_{\rho\times \rho}$. This proves \eqref{st:nowall}.
\medskip

\sta{\label{st:nocompbip} $G$ has no induced minor isomorphic to $K_{\sigma,\sigma}$.}

By \eqref{st:nowall}, $G$ has no induced minor isomorphic to $W_{\rho\times \rho}$. Also, note that there is a theta in every $(2,3)$-constellation. Since $G$ is theta-free, it follows that there is no $(2,3)$-constellation in $G$. Thus, by \Cref{thm:compbip} and the choice of $\sigma$, $G$ has no induced minor isomorphic to $K_{\sigma,\sigma}$. This proves \eqref{st:nocompbip}.
\medskip

Note that if $t=1$, then $\tw(G)<1=t^{d_{\ref{thm:main}}}$, as desired. Assume that $t\geq 2$. Then, since $G$ is $t^d$-separable, from \eqref{st:nowall}, \eqref{st:nocompbip}, \Cref{thm:mainpolyblock}, and the choice of $d'$, it follows that 
$$\tw(G)\leq \left(2\cdot t^{d}\cdot t\right)^{d'}\leq \left(t^{(d+2)}\right)^{d'}=t^{d_{\ref{thm:main}}}.$$
This completes the proof of \Cref{thm:main}.
\end{proof}

\section{Growing a tree}\label{sec:final}

In this section, we prove \Cref{thm:mainseparation}. For the purpose of an inductive argument, we need to prove the following slightly stronger lemma. For $a,b\in \poi$, by an \textit{$(a,b)$-tree} we mean a pair $(T,x)$ where $T$ is a tree and $x\in V(T)$ such that, if $b=1$, then $V(T)=\{x\}$, and if $b\geq 2$, then $x$ has degree exactly $a$ in $T$, every non-leaf vertex of $T$ other than $x$ has degree exactly $a$ in $T$, and every leaf of $T$ is at distance exactly $b-1$ from $x$ in $T$.  Given a graph $G$ and a set $\mca{P}$ of paths in $G$, we write $V(\mca{P})=\bigcup_{P\in \mca{P}}V(P)$. 

\begin{lemma}\label{lem:treeramsey}
For all $a,b\in \poi$, there are constants $c_{\ref{lem:treeramsey}}=c_{\ref{lem:treeramsey}}(a,b)\in \poi$ and $d_{\ref{lem:treeramsey}}=d_{\ref{lem:treeramsey}}(a,b)$ with the following property. Let $t\in \poi$ and let $G$ be a theta-free $K_{t}$-free graph. Let $x,y\in V(G)$ be distinct and nonadjacent, and let $\mca{P}$ be a set of pairwise internally disjoint paths in $G$ from $x$ to $y$ such that $|\mca{P}|\geq c_{\ref{lem:treeramsey}}t^{d_{\ref{lem:treeramsey}}}$. Then $G[V(\mca{P})\setminus \{y\}]$ has an induced subgraph $T$ with $x\in V(T)$ such that $(T,x)$ is an $(a,b)$-tree.
\end{lemma}

The proof uses the following from \cite{ti6} (the exponent $8s+12$ is not explicitly given in the statement of Theorem~4.2 in \cite{ti6}, but it may easily be recovered from the proof). A set $Z$ of vertices in a graph $G$ is \textit{constricted} if $Z$ is a stable set in $G$ with $|Z|\geq 3$ and no three vertices in $Z$ belong to an induced subgraph of $G$ that is a tree.

\begin{theorem}[Chudnovsky and Codsi; Theorem~4.2 in \cite{ti6}]\label{thm:3inatreesep}
    Let $s\in \poi$, let $G$ be a $K_{s,s}$-free graph, let $Z\subseteq V(G)$ be constricted in $G$, and let $y\in V(G)\setminus Z$. Then there are at most $\omega(G)^{8s+12}$ pairwise disjoint paths (except at $y$) from $y$ to a vertex in $Z$.
\end{theorem}

In addition, we need several Ramsey-type results. The first is about digraphs. By a \textit{digraph} we mean a pair $D=(V(D), E(D))$ where $V(D)$ is a finite set of \textit{vertices} and $E(D)\subseteq (V(D)\times V(D))\setminus \{(v,v):v\in V(D)\}$ is the set of \textit{edges}. In particular, our digraphs are loopless and allow at most one edge in each direction between every two vertices. For $(u,v)\in E(D)$, we say that $v$ is an \textit{out-neighbor} of $u$, and the \textit{out-degree} of a vertex $u\in V(D)$ is the number of its out-neighbors. The \textit{underlying graph of $D$} is the graph $G$ with $V(G)=V(D)$ and $E(G)=\{uv:(u,v)\in E(D) \text{ or } (v,u)\in E(D)\}$. A \textit{stable set in $D$} is a stable set in the underlying graph of $D$.

\begin{lemma}[Chudnovsky, Hajebi, Spirkl; Lemma 5.1 in \cite{tw18}]\label{lem:digraph}
    Let $q,r,s\in \poi$ and let $D$ be a digraph. Then the following hold.
    \begin{enumerate}[{\rm (a)}]
    \item \label{lem:digraph_a} 
    If $D$ has at least $2rs$ vertices of out-degree at most $r$, then there is a stable set of cardinality $s$ in $D$.
    \item \label{lem:digraph_b} If there are at least $2qrs$ vertices of out-degree at least $qr$ in $D$, then there is an $s$-subset $S$ of $V(D)$ with the following property: for every $q$-subset $\{v_1,\ldots, v_q\}$ of $S$, there are $q$ pairwise disjoint $r$-subsets $R_1,\ldots, R_q$ of $V(D)\setminus S$ such that for each $i\in \{1,\ldots, q\}$, every vertex in $R_i$ is an out-neighbor of $v_i$.
    \end{enumerate}
\end{lemma}

We also use the classical theorem of Ramsey \cite{multiramsey}, as well as an improvement of it when restricted to graphs that exclude a fixed complete bipartite graph. (For the familiar reader, \Cref{lem:EHforKss} simply restates the well-known fact that the Erd\H{o}s-Hajnal conjecture \cite{EH} is true for complete bipartite graphs, which, for instance, follows from the main result of \cite{EH2}. See Lemma~5.17 in \cite{polypw} for a short proof with the explicit bound.)

\begin{theorem}[Ramsey \cite{multiramsey}]\label{thm:ramsey}
    Let $\alpha,t\in \poi$ and let $G$ be a $K_t$-free graph with no stable set of cardinality $\alpha$. Then $|V(G)|<t^{\alpha-1}$.
\end{theorem}

\begin{lemma}[Erd\H{o}s and Hajnal \cite{EH2}; See also Lemma~5.17 in \cite{polypw}]\label{lem:EHforKss}
    Let $\alpha, s, t\in \poi$ and let $G$ be a $(K_{s,s},K_{t})$-free graph with no stable set of cardinality $\alpha$. Then $|V(G)|< \alpha^st^{s-1}$.
\end{lemma}

Next, we prove a strengthening of \Cref{lem:EHforKss}. (A similar result with worse bounds is also proved in \cite{lozin}.)
For a graph $G$, we say that $X,Y\subseteq V(G)$ \textit{are anticomplete in $G$} if $X\cap Y=\varnothing$ and there is no edge of $G$ with an end in $X$ and an end in $Y$. (Note that \Cref{lem:EHforKss} is identical to \Cref{lem:smallsetsanti} below for $r=1$.) 

\begin{lemma}\label{lem:smallsetsanti}
 Let $r,s\in \poi$ and let $\sigma_{r}=s^{(4s)^{r-1}}$. Let $\alpha,t\in \poi$, let $G$ be a $(K_{s,s},K_t)$-free graph and let $\mca{X}$ be a set of pairwise disjoint nonempty subsets of $V(G)$, each of cardinality at most $r$, such that  $|\mca{X}|\geq \alpha^{\sigma_{r}}t^{\sigma_{r}-1}$. Then there are $\alpha$ elements in $\mca{X}$ that are pairwise anticomplete in $G$.
\end{lemma}
\begin{proof}
The result is trivial for $\alpha=1$. Also, if $s=1$, then $|\mca{X}|\geq \alpha$ and the sets in $\mca{X}$ are all pairwise anticomplete in $G$, as desired. Assume that $\alpha,s\geq 2$. Let us now show that:

\sta{\label{st:sigmabounds} We have
$$\alpha^{\sigma_{2}}t^{\sigma_{2}-1}> \alpha^{s}t^{s-1}+\alpha^{\left((2s)^{2s-1}-1\right)s^2}t^{s^{2}-1}.$$
Moreover, for all $r\in \poi$ with $r\geq 2$, we have
$$\alpha^{\sigma_{r}}t^{\sigma_{r}-1}> \alpha^{\sigma_{r-1}}t^{\sigma_{r-1}-1}+\alpha^{\sigma^3_{r-1}}t^{\sigma_{r-1}^3-1}.$$}

Since $s\geq 2$, it follows that $\sigma_{2}=s^{4s}=s^2+((s^2)^{2s-1}-1)s^2\geq 1+((2s)^{2s-1}-1)s^2$. Thus, since $\alpha\geq 2$, it follows that 
   $$\alpha^{\sigma_{2}}> \alpha^{1+\left((2s)^{2s-1}-1\right)s^2}\geq 2\alpha^{\left((2s)^{2s-1}-1\right)s^2}\geq \alpha^s+ \alpha^{\left((2s)^{2s-1}-1\right)s^2}.$$
In particular, since $\sigma_{2}>s^2>s\geq 2$, it follows that
$$\alpha^{\sigma_{2}}t^{\sigma_{2}-1}> \alpha^{s}t^{s-1}+\alpha^{\left((2s)^{2s-1}-1\right)s^2}t^{s^{2}-1}.$$
This proves the first assertion. Let $r\geq 2$. Since $s\geq 2$, it follows that $\sigma_{r-1}=s^{(4s)^{r-2}}\geq 2$ and so
$$\sigma_{r}=\sigma_{r-1}^{4s}>\sigma_{r-1}^4>\sigma_{r-1}^3+1.$$
Therefore, since $\alpha\geq 2$, it follows that:
$$\alpha^{\sigma_{r}}t^{\sigma_{r}-1}>\alpha^{\sigma^3_{r-1}+1}t^{\sigma_{r-1}^3}\geq 2\alpha^{\sigma^3_{r-1}}t^{\sigma_{r-1}^3}\geq \alpha^{\sigma_{r-1}}t^{\sigma_{r-1}-1}+\alpha^{\sigma^3_{r-1}}t^{\sigma_{r-1}^3-1}.$$
This proves \eqref{st:sigmabounds}.
\medskip

Now, we prove \ref{lem:smallsetsanti} by induction on $r$, for fixed $s,t\in \poi$ with $s\geq 2$, and all $\alpha\in \poi$ with $\alpha\geq 2$. If $r=1$, then $\sigma_{r}=s$, $|\mca{X}|\geq \alpha^{s}t^{s-1}$, and the result follows from \Cref{lem:EHforKss}. So we may assume that $r\geq 2$. Let $\mca{X}^-\subseteq \mca{X}$ be the set of all sets in $\mca{X}$ of cardinality at most $r-1$. Note that if $|\mca{X}^-|\geq \alpha^{\sigma_{r-1}}t^{\sigma_{r-1}-1}$, then we are done by applying the inductive hypothesis to $\mca{X}^-$. Therefore:

\sta{\label{st:smallxminus}  We may assume that $|\mca{X}^-|< \alpha^{\sigma_{r-1}}t^{\sigma_{r-1}-1}$.}

We further prove that:

\sta{\label{st:r=2}  We may assume that $r\geq 3$.}

Assume that $r=2$. Let $\zeta_2=\alpha^{(2s)^{2s-1}-1}$, let $\zeta_1=\zeta_2^{s}t^{s-1}$ and let $\zeta_0=\zeta_1^{s}t^{s-1}$. Then
$$\zeta_0=\left(\zeta_2^{s}t^{s-1}\right)^st^{s-1}=\left(\left(\alpha^{(2s)^{2s-1}-1}\right)^st^{s-1}\right)^st^{s-1}=\alpha^{\left((2s)^{2s-1}-1\right)s^2}t^{s^{2}-1}.$$
This, along with \eqref{st:sigmabounds} and \eqref{st:smallxminus}, yields 
$|\mca{X}\setminus \mca{X}^-|>\alpha^{\sigma_{2}}t^{\sigma_{2}-1}-\alpha^{s}t^{s-1}>\alpha^{\left((2s)^{2s-1}-1\right)s^2}t^{s^{2}-1}=\zeta_0.$
Note that every set in $\mca{X}\setminus \mca{X}^-$ has cardinality exactly $r=2$. Choose $\zeta_0$ pairwise distinct sets $\{x_1,y_1\},\ldots, \{x_{\zeta_0},y_{\zeta_0}\}\in \mca{X}\setminus \mca{X}^-$. Since $G$ is $(K_{s,s},K_t)$-free and $\zeta_0=\zeta_1^st^{s-1}$, by \Cref{lem:EHforKss} applied to $G[\{x_1,\ldots, x_{\zeta_0}\}]$, there exists $I_1\subseteq \{1,\ldots, \zeta_0\}$ with $|I_1|=\zeta_1$ such that $\{x_i:i\in I_1\}$ is a stable set in $G$. Similarly, since $|I_1|=\zeta_1=\zeta_2^st^{s-1}$, by \Cref{lem:EHforKss} applied to $G[\{y_i:i\in I_1\}]$, there exists $I_2\subseteq I_1$ with $|I_2|=\zeta_2$ such that $\{y_i:i\in I_2\}$ is a stable set in $G$. Let $\Gamma$ be the graph with $V(\Gamma)=I_2$ such that for all $i,j\in I_2$ with $i<j$, we have $ij\in E(\Gamma)$ if and only if $x_iy_j,x_jy_i\notin E(G)$. We claim that $\Gamma$ has a clique of cardinality $\alpha$. Suppose not. Then, since $|V(\Gamma)|=|I_2|=\zeta_2=\alpha^{(2s)^{2s-1}-1}$, by \Cref{thm:ramsey}, there is a stable set in $\Gamma$ of cardinality $(2s)^{2s-1}$; that is, there exists $I\subseteq I_2\subseteq I_1$ with $|I|=(2s)^{2s-1}$ such that for all $i,j\in I_2$ with $i<j$, either $x_iy_j\in E(G)$ or $x_jy_i\in E(G)$. Let $\Gamma'$ be the graph with $V(\Gamma')=I$ such that for all $i,j\in I$ with $i<j$, we have $ij\in E(\Gamma')$ if and only if $x_iy_j\in E(G)$. Since $|V(\Gamma')|=|I|=(2s)^{2s-1}$, by \Cref{thm:ramsey}, there is either a clique or a stable set of cardinality $2s$ in $\Gamma'$; in particular, there are $i_1,\ldots, i_{s}, j_1,\ldots, j_s\in  I\subseteq I_2\subseteq I_1$ with $i_1<\cdots<i_{s}< j_1<\cdots<j_s$ such that, either for all $k,l\in \{1,\ldots, s\}$, we have $x_{i_k}y_{j_l}\in E(G)$, or for all $k,l\in \{1,\ldots, s\}$, we have $x_{j_k}y_{i_l}\in E(G)$. But in the former case $G[\{x_{i_1},\ldots, x_{i_{s}}, y_{j_1},\ldots, y_{j_s}\}]$ is isomorphic to $K_{s,s}$, and in the latter case, $G[\{x_{j_1},\ldots, x_{j_{s}}, y_{i_1},\ldots, y_{i_s}\}]$ is isomorphic to $K_{s,s}$, contrary to the assumption that $G$ is $K_{s,s}$-free. This proves the claim that $\Gamma$ has a clique of cardinality $\alpha$; that is, there exists $I_3\subseteq I_2\subseteq I_1$ with $|I|=\alpha$ such that for all $i,j\in I_2$ with $i<j$, we have $x_iy_j, x_jy_i\notin E(G)$. But now $(\{x_i,y_i\}:i\in I_3)$ are $\alpha$ sets in $\mca{X}\setminus \mca{X}^-$ that are pairwise anticomplete in $G$, as required. This proves \eqref{st:r=2}.
\medskip

 Let $\xi_2=\alpha^{\sigma_{r-1}}t^{\sigma_{r-1}-1}$, let $\xi_1=\xi_2^{\sigma_{r-1}}t^{\sigma_{r-1}-1}$ and let $\xi_0=\xi_1^{\sigma_{r-1}}t^{\sigma_{r-1}-1}$. Then
$$\xi_0=\left(\xi_2^{\sigma_{r-1}}t^{\sigma_{r-1}-1}\right)^{\sigma_{r-1}}t^{\sigma_{r-1}-1}=\left(\left(\alpha^{\sigma_{r-1}}t^{\sigma_{r-1}-1}\right)^{\sigma_{r-1}}t^{\sigma_{r-1}-1}\right)^{\sigma_{r-1}}t^{\sigma_{r-1}-1}=\alpha^{\sigma^3_{r-1}}t^{\sigma_{r-1}^3-1}.$$
Thus, by \eqref{st:sigmabounds} and \eqref{st:smallxminus}, we have
$|\mca{X}\setminus \mca{X}^-|>\alpha^{\sigma_{r}}t^{\sigma_{r}-1}-\alpha^{\sigma_{r-1}}t^{\sigma_{r-1}-1}>\alpha^{\sigma^3_{r-1}}t^{\sigma_{r-1}^3-1}=\xi_0$. Choose $\xi_0$ pairwise distinct sets $W_1,\ldots, W_{\xi_0}\in \mca{X}\setminus \mca{X}^-$. Then, $$|W_1|=\cdots=|W_{\xi_0}|=r.$$
 
 By \eqref{st:r=2}, we may assume that $r\geq 3$. For each $i\in \{1,\ldots, \xi_0\}$, choose three pairwise distinct vertices $x_i,y_i,z_i$ in $W_i$, let $X_i=W_i\setminus \{x_i\}$, let $Y_i=W_i\setminus \{y_i\}$ and let $Z_i=W_i\setminus \{z_i\}$. Since $G$ is $(K_{s,s},K_t)$-free, $|X_1|=\cdots=|X_{\xi_0}|=r-1$ and $\xi_0=\xi_1^{\sigma_{r-1}}t^{\sigma_{r-1}-1}$, by the inductive hypothesis applied to $\mca{X}_0=\{X_1,\ldots, X_{\xi_0}\}$, there exists $I_1\subseteq \{1,\ldots, \xi_0\}$ with $|I_1|=\xi_1$ such that $(X_i:i\in I_1)$ are pairwise anticomplete in $G$. Similarly, since $|Y_i|=r-1$ for every $i\in I_1$, and since $|I_1|=\xi_1=\xi_2^{\sigma_{r-1}}t^{\sigma_{r-1}-1}$, by the inductive hypothesis applied to $\mca{X}_1=\{Y_i:i\in I_1\}$, there exists $I_2\subseteq I_1$ with $|I_2|=\xi_2$ such that $(Y_i:i\in I_2)$ are pairwise anticomplete in $G$. Finally, since $|Z_i|=r-1$ for every $i\in I_2$, and since $|I_2|=\xi_2=\alpha^{\sigma_{r-1}}t^{\sigma_{r-1}-1}$, by the inductive hypothesis applied to $\mca{X}_2=\{Z_i:i\in I_2\}$, there exists $I_3\subseteq I_2\subseteq I_1$ with $|I_3|=\alpha$ such that $(Z_i:i\in I_3)$ are pairwise anticomplete in $G$. Hence, $(W_i=X_i\cup Y_i\cup Z_i:i\in I_3)$ are $\alpha$ sets in $\mca{X}\setminus \mca{X}^-$ that are pairwise anticomplete in $G$. This completes the proof of \Cref{lem:smallsetsanti}.
\end{proof}

Now we can prove \Cref{lem:treeramsey}:

\begin{proof}[Proof of \Cref{lem:treeramsey}]

Throughout, let $a\in \poi$ be fixed. For every $n\in \poi$, let $$\theta_n=3^{12^{na^{n-1}-1}}.$$
We define the sequences $(\mu_{n}:n\in \poi)$ and $(\lambda_{n}:n\in \poi)$ recursively, as follows: Let $\mu_1=\lambda_1=0$, and for $n\in \poi$ with $n\geq 2$, let
$$\mu_n=\left(\left(3 a^{2\theta_n}+2\right) \mu_{n-1}\right)^3$$
and let
$$\lambda_n=3\lambda_{n-1}+6\theta_n-4.$$

We will prove by induction on $b\in \poi$ that 
$$c_{\ref{lem:treeramsey}}=c_{\ref{lem:treeramsey}}(a,b)=\mu_b, \quad d_{\ref{lem:treeramsey}}=d_{\ref{lem:treeramsey}}(a,b)=\lambda_b$$
satisfy the lemma.

Let $(T,v)$ be an $(a,b)$-tree. Let $t\in \poi$ and let $G$ be a theta-free $K_{t}$-free graph. Since $G$ is theta-free, it follows that $G$ is $K_{3,3}$-free. Let $x,y\in V(G)$ be distinct and nonadjacent, and let $\mca{P}$ be a set of pairwise internally disjoint paths in $G$ from $x$ to $y$ such that $|\mca{P}|\geq \mu_b t^{\lambda_b}\geq 1$. Note that if $b=1$, then $T=G[\{x\}]$ is the desired induced subgraph of $G[V(\mca{P})\setminus \{y\}]$. Moreover, since $\mca{P}\neq \varnothing$, it follows that $G$ is not edgeless. So, we may assume that $b\geq 2$ and $t\geq 3$. 

For each $P\in \mca{P}$, let $x_P$ be the unique neighbor of $x$ in $P$ (so $x_P\neq y$). For every subset $\mca{Q}$ of $\mca{P}$, let $X_{\mca{Q}}=\{x_Q:Q\in \mca{Q}\}$. Note that
$$|\mca{P}|\geq \mu_b t^{\lambda_b}= \left(\left(3 a^{2\theta_b}+2\right) \mu_{b-1}\right)^3\cdot t^{3\lambda_{b-1}+6\theta_b-4}=\left(\left(3 a^{2\theta_b}+2\right)\cdot \left(t^{2(\theta_b-1)}\mu_{b-1} t^{\lambda_{b-1}}\right)\right)^3t^2.$$

Since $G$ is $(K_{3,3}, K_{t})$-free (because $G$ is theta-free and $K_{t}$-free), by \Cref{lem:EHforKss} applied to $G[X_{\mca{P}}]$, there exists 
$\mca{Q}\subseteq \mca{P}$ with 
$$
|\mca{Q}|=\left(3 a^{2\theta_b}+2\right)\cdot \left(t^{2(\theta_b-1)}\mu_{b-1} t^{\lambda_{b-1}}\right)\geq 3\left(a^{\theta_b} t^{\theta_b-1}\right)^2\cdot \mu_{b-1} t^{\lambda_{b-1}}+2$$
such that $X_{\mca{Q}}$ is a stable set in $G$. In particular, if there are three paths $Q_1,Q_2,Q_3\in \mca{Q}$ of length two, then $G[Q_1\cup Q_2\cup Q_3]$ is a theta in $G$, a contradiction. It follows that there exists $\mca{L}\subseteq \mca{Q}\subseteq \mca{P}$ with 
$$|\mca{L}|=3\left(a^{\theta_b} t^{\theta_b-1}\right)^2\cdot \mu_{b-1} t^{\lambda_{b-1}}$$
such that $X_{\mca{L}}$ is a stable set in $G$ and every path in $\mca{L}$ has length at least three in $G$; that is, $N_{X_{\mca{L}}}(y)=\varnothing$, and so for every $P\in \mca{L}$, the vertex $x_P$ has a unique neighbor $x'_P$ in $P^*$. 

Let 
$q=a^{\theta_{b}} t^{\theta_{b}-1}$ and let
$r=\mu_{b-1} t^{\lambda_{b-1}}$. Since $b\geq 2$, it follows that $q\geq 2^{3^{12}}t^{3^{12}-1}> 2t^{36}$, and so
$|\mca{L}|\geq 3q^2r\geq (2q+2t^{36})qr$. Let $D$ be the digraph with $V(D)=\mca{L}$ such that for distinct $P,Q\in \mca{L}$, we have $(P,Q)\in E(D)$ if and only if $N_{Q^*}(x_{P})\neq \varnothing$.

  \sta{\label{st:nostableset} There are at least $2q^2r$ vertices in $D$ of out-degree at least $qr$.}
  
  Suppose not. Then, since $|V(D)|=|\mca{L}|\geq (2q+2t^{36})qr$, it follows that there are at least $2t^{36}qr$ vertices in $D$ with out-degree at most $qr$. Thus, by \Cref{lem:digraph}\ref{lem:digraph_a}, there is a stable set $\mca{S}\subseteq \mca{L}=V(D)$ in $D$ with $|\mca{S}|=t^{36}$.
  Let $G_1=G[V(\mca{S})\setminus \{x\}]$. By the definition of $D$, for every $P \in \mca{S}$, we have $N_{G_1}(x_{P})=N_{P^*}(x_P)=\{x'_P\}$, and, in particular, $|N_{G_1}(x_{P})|=1$. Let $Z=X_{\mca{S}}$. Then $Z$ is a stable set in $G_1$, $y\in G_1\setminus Z$, and $\{P\setminus \{x\}:P\in \mca{S}\}$ is a set of $|Z|=|\mca{S}|=t^{36}>\omega(G)^{36}$ paths of nonzero length in $G$, all sharing $y$ as an end and otherwise pairwise disjoint, and each having an end in $Z$. Since $G_1$ is $(K_{3,3},K_{t})$-free and $|Z|\geq 3$ (because $t\geq 2$), it follows from \Cref{thm:3inatreesep} that $Z$ is not constricted in $G_1$. Consequently, there is an induced subgraph of $G_1$ which is a tree that intersects $Z$ in three vertices. In particular, there is an induced subgraph $T$ of $G_1$ which is a tree with $|T\cap Z|=3$, say $T\cap Z=\{x_{P_1},x_{P_2},x_{P_3}\}$ where $P_1,P_2,P_3\in \mca{S}$. Since $x_{P_1},x_{P_2},x_{P_3}$ all have degree one in $G_1$ (and so in $T$), it follows that that there is an induced subgraph $T'$ of $T$ isomorphic to a subdivision $K_{1,3}$ and $x_{P_1},x_{P_2},x_{P_3}$ are the degree-one vertices of $T'$. Note also that $V(T')\setminus \{x_{P_1},x_{P_2},x_{P_3}\}\subseteq V(G_1)\setminus Z$, and so $x$ is anticomplete to $V(T')\setminus \{x_{P_1},x_{P_2},x_{P_3}\}$ in $G$. But now 
$G[V(T')\cup \{x\}]$
is a theta in $G$, a contradiction. This proves
\eqref{st:nostableset}.
  \medskip

By \eqref{st:nostableset},  \Cref{lem:digraph}\ref{lem:digraph_b} and the definition of $D$, there are $q$ pairwise disjoint $(r+1)$-subsets $\mca{L}_1, \ldots, \mca{L}_{q}$ of $\mca{L}$, such that for each $i\in \{1,\ldots, q\}$, one may write $\mca{L}_i=\{P_i\}\cup \mca{R}_i$
with $|\mca{R}_i|=r$, such that $x_{P_i}$ has a neighbor in $R^*$ for every $R\in \mca{R}_i$. Thus, for each $i\in \{1,\ldots, q\}$ and every $R\in \mca{R}_i$, there is a path $P_R$ in $G$ from $x_{P_i}$ to $y$ with $P^*_R\subseteq R^*$.

For each $i\in \{1,\ldots, q\}$, let
$\mca{P}_i=\{P_R:R\in \mca{R}_i\}$. Then $\mca{P}_i$ is a set of $r=\mu_{b-1}t^{\lambda_{b-1}}$ pairwise internally disjoint paths in $G$ from $x_{P_i}$ to $y$. For every $i\in \{1,\ldots, q\}$, since $x_{P_i},y$ are nonadjacent in $G$, it follows from the inductive hypothesis that $G[V(\mca{P}_i)\setminus \{y\}]$ has an induced subgraph $T_i$ with $x_{P_i}\in V(T_i)$ such that $(T_i,x_{P_i})$ is an $(a,b-1)$-tree. In particular, $V(T_1),\ldots,V(T_{q})$ are pairwise disjoint (because $V(\mca{P}_1)\setminus \{y\}, \ldots, V(\mca{P}_{q})\setminus \{y\}$ are pairwise disjoint), and for each $i\in \{1,\ldots, q\}$, we have $N_{V(T_i)}(x)=N_{V(\mca{P}_i)}(x)=\{x_{P_i}\}$. 

Let $r=ba^{b-1}$. Then $$q=a^{\theta_b} t^{\theta_b-1}=a^{3^{12^{r-1}}}t^{3^{12^{r-1}}-1}.$$
Therefore, since $G$ is $(K_{3,3}, K_t)$-free, since $|V(T_i)|=|V(T')|=\sum_{i=0}^{b-1}a^i\leq r$ for all $i\in \{1,\ldots, q\}$, it follows from \Cref{lem:smallsetsanti} applied to $\mca{X}=\{V(T_1),\ldots, V(T_q)\}$ that there is an $a$-subset $I$ of $\{1,\ldots, q\}$ for which $(V(T_i):i\in I)$ are pairwise anticomplete in $G$. 

Now, let $$T=G\left[\left(\bigcup_{i\in I}V(T_{i})\right)\cup \{x\}\right].$$
Then $T$ is an induced subgraph of $G[V(\mca{P})\setminus \{y\}]$ with $x\in V(T)$ such that $(T,x)$ is an $(a,b)$-tree. This completes the proof of \Cref{lem:treeramsey}.
\end{proof}

We are now in a position to prove \Cref{thm:mainseparation}, which we restate:

\mainseparation*

\begin{proof}

Let $c=c_{\ref{lem:treeramsey}}(h+1,h+1)$, let $d=d_{\ref{lem:treeramsey}}(h+1,h+1)$, and let
$$d_{\ref{thm:mainseparation}}=d_{\ref{thm:mainseparation}}(h)=c+d.$$

Let $G$ be a theta-free $(H, K_{t})$-free graph that is not $t^{d_{\ref{thm:mainseparation}}}$-separable. Then there are distinct and nonadjacent vertices $x,y\in V(G)$ for which there is a set $\mca{P}$ of pairwise internally disjoint paths in $G$ from $x$ to $y$ such that $|\mca{P}|\geq t^{d_{\ref{thm:mainseparation}}}\geq 1$. In particular, since $\mca{P}\neq \varnothing$, it follows that $t\geq 3$, and so $|\mca{P}|\geq t^{d_{\ref{thm:mainseparation}}}=t^{c+d}\geq 3^{c}t^{d}> c t^{d}$. By \Cref{lem:treeramsey}, $G$ has an induced subgraph $T$ with $x\in V(T)$ such that $(T,x)$ is an $(h+1,h+1)$-tree. Moreover, since $H$ is a forest on at most $h$ vertices, it follows that there is a tree $H^+$ on at most $h+1$ vertices such that $H$ is an induced subgraph of $H^+$, and since $|V(H^+)|\leq h+1$, it follows that $T$ has an induced subgraph isomorphic to $H^+$. But now $G$ has an induced subgraph isomorphic to $H$, a contradiction. This completes the proof of \Cref{thm:mainseparation}.
\end{proof}

\bibliographystyle{plain}
\bibliography{ref}

\end{document}